\documentclass{ifacconf}
\usepackage{graphicx}      % include this line if your document contains figures
\usepackage{natbib}        % required for bibliography
\usepackage{amsmath,amsfonts}
\def\x{\mathbf{x}}

\def\P{\mathbf{P}}
\def\R{\mathbb{R}}
\def\bl{\boldsymbol{\ell}}

\def\x{\mathbf{x}}
\def\P{\mathbf{P}}
\def\S{\mathbf{S}}

\def\y{\mathbf{y}}

\def\R{\mathbb{R}}

\def\bH{\mathbf{H}}
\def\N{\mathbb{N}}

\def\Q{\mathbf{Q}}
\def\M{\mathbf{M}}

\def\E{\mathcal{E}}

\def\A{\mathbf{A}}

\def\L{\mathcal{L}}
\def\bS{\mathbb{S}}

\def\z{\mathbf{z}}
\def\v{\mathbf{v}}
\def\y{\mathbf{y}}

\def\x{\mathbf{x}}
\def\bs{\mathbf{s}}

\def\y{\mathbf{y}}

\def\b{\mathbf{b}}

\def\u{\mathbf{u}}

\def\e{\mathbf{e}}

\def\bom{\boldsymbol{\Omega}}
\def\bla{\boldsymbol{\lambda}}
\def\l{\boldsymbol{\ell}}

\begin{document}
\begin{frontmatter}

\title{Optimization of Polynomials with Sparsity Encoded in a Few Linear Forms\thanksref{footnoteinfo}} 

\thanks[footnoteinfo]{Research sponsored by the Artificial and Natural Intelligence Institute (ANITI) of Toulouse,
and ANR-NuSCAP-20-CE48-0014}

\author[First]{Jean B. Lasserre} 

\address[First]{LAAS-CNRS \& Institute of Mathematics, University of Toulouse, France (e-mail: lasserrer@ laas.fr).}

\begin{abstract}   
We consider polynomials of a few linear forms and show how exploit this 
 type of sparsity for optimization on some particular domains like the Euclidean sphere or a polytope.
 Moreover, a simple procedure allows to detect this form of sparsity and also allows to provide an approximation of any polynomial by such sparse polynomials.
 \end{abstract}


\end{frontmatter}
%===============================================================================

\section{Introduction}
In this paper we discuss the optimization problems 
\[\P:\quad\min \{h(\x): \x\in\bom\}\,,\]
where 
the polynomial\footnote{Most of what follows also applies to continuously 
differentiable functions} $h\in\R[x_1,\ldots,x_n]$ is a 
defined in terms of a few linear forms, that is,
\begin{equation}
\label{def-sparse}
\x\mapsto h(\x)\,=\,f(\bl^T\x)\,,\quad \x\in\R^n\,,\end{equation}
for some polynomial $f\in\R[X_1,\ldots,X_m]$ 
and some real matrix $\bl\in\R^{n\times m}$. 
We are also interested in approximating an arbitrary polynomial 
by polynomials in the form \eqref{def-sparse}.

\subsubsection{Motivation.} When $m\ll n$, formulation \eqref{def-sparse} exhibits some sort of \emph{sparsity} as only a few linear forms are involved in $h$. Indeed such a sparsity has been explored in several contexts
like e.g. statistical learning in \cite{Roweis}, \cite{camastra},
to learn a low-dimensional manifold 
(where $h$ is called a low-rank function), 
in \cite{active-set} for contracting response surfaces
on a low-dimensional subspace, in \cite{vergne} for multivariate
integration on the simplex, and in \cite{barvinok} for integration with respect to the Gaussian measure. Therefore 
one also expects that it can be exploited for an efficient computation of 
the (local or global) minimum on $\bom$. But notice that in general
$\bom$ is \emph{not} expressed in terms of the $\bl_j$'s so that exploiting this sparsity
to optimize $h$ on $\bom$ may not be easy.  For the sphere $\bS^{n-1}$ 
\cite{barvinok} has shown that the maximum (but not the minimum) 
of certain sparse homogeneous polynomials 
can be approximated well by a properly scaled maximum on the unit sphere of a
random low-dimensional subspace. This class of homogeneous polynomials contains some polynomials of the form \eqref{def-sparse}. Notice also that the sparsity \eqref{def-sparse} (when $m\ll n$)  is different from the various sparsity patterns exploited for polynomial optimization in \cite{amirali}, \cite{lass-2006}, \cite{tssos1}, and \cite{tssos2}.

If $h$ is not directly available in sparse form \eqref{def-sparse}, its detection
is quite important in view of the potential resulting benefits for optimization.  It turns out that the detection issue has been already addressed in engineering and data science, in the more general context of approximating
an arbitrary continuous differentiable function 
$h(\bl \y+\bs\,\v)$ where the columns of $\bl\in\R^{n\times m}$,
(resp. $\bs\in\R^{n\times (n-m)}$) are eigenvectors of $\mathrm{E}_\mu [\nabla h\nabla h^T]$ associated with the 
$m$ largest (resp. the remaining $n-m$) eigenvalues, and $\mu$
is an appropriate probability measure. 
%The idea is to approximate $h(\x)$ with $G(\bl^T\x)$, where $G(\y):=\mathrm{E}_\mu[\,h\vert\y\,]$.
In \cite{active-set} the authors 
discuss methods to obtain and evaluate an approximation
based on the function $G(\y):=\mathrm{E}_\mu[\,h\vert\y\,]$; see below. (In \cite{active-set} $h$ in \eqref{def-sparse} is called a $\z$-invariant function.) 
Notice that if even if $h$ is a polynomial,
the resulting approximation $G$ is \emph{not}.

\subsubsection{Contribution.} Our contribution is threefold:

(i) We show that the sparsity in \eqref{def-sparse} can be exploited in optimization on the Euclidean sphere $\bS^{n-1}$
%, the Euclidean ball $\E_n$, 
and arbitrary polytopes.
%, and the discrete hypercube $\{-1,1\}^n$.
Solving the original problem reduces to solving 
an explicit optimization problem in $\R^m$,
simply related and similar to $\P$, but with a drastic reduction in difficulty. We thus extend \cite{lass-sphere} who considered the case $\bom=\mathbb{S}^{n-1}$
and showed that solving $\P$ is equivalent to minimizing the $m$-variables polynomial
$X\mapsto f(\L_1\cdot X_1,\ldots,\L_m\cdot X_m)$
on the Euclidean ball $\mathcal{E}_m$,
(where $\L_i$ is the $i$-th column of $\bl^T\bl$). 

(ii) A second contribution is with respect to \emph{detection} of 
a sparsity \eqref{def-sparse}. When $h$ is a polynomial we provide two procedures. We first choose $\mu$ to be the uniform distribution on $\E_n$. Then we build a matrix
$\bH_k\bH_k^T$ where the columns of $\bH_k$ are just the gradient of $h$ evaluated at points $(\x(1),\ldots,\x(k))\subset\E_n$ (randomly generated according to $\mu$), until the condition ${\rm rank}(\bH_k\bH_k^T)=
{\rm rank}(\bH_{k-1}\bH_{k-1}^T)$ is satisfied, say for $k=m+1$. Then a  sparsity as in \eqref{def-sparse} for some explicit $\bl\in\R^{n\times m}$, is detected with probability $1$. A second possibility that gets rid of ``\emph{with prob.} $1$" is 
to follow \cite{active-set} and perform the SVD decomposition of 
$\mathcal{R}:=\mathrm{E}_\mu[\nabla h\,\nabla h^T]$. 
But in pour context, as $h$ is a polynomial and integration of polynomials on $\E_n$ is easy, $\mathcal{R}$ can be computed \emph{exactly}. Then $\bl$ in \eqref{def-sparse} is obtained from eigenvectors associated with the $m$ non-zero eigenvalues of $\mathcal{R}$.

(iii) A third contribution is with respect to detection of an \emph{approximate sparsity} and is directly inspired by the active set method, as described in e.g. \cite{active-set}. Write $h$ in the form
\[h(\x)\,=\,h(\bl\,\y+\bs\,\z)\,,\quad\bl\in\R^{n\times m}\,\:\bs\in\R^{n\times (n-m)}\,,\]
where the columns of $\bl,\bs$ are the eigenvectors of  
$\mathcal{R}$ (with norm $1$), and where the $n-m$ eigenvalues associated with $\bs$
are much smaller than the $m$ eigenvalues associated with $\bl$.
 In \cite{active-set} the authors propose to approximate 
$h$ with the function $G(\bl^T\x)$ defined by:
\begin{equation}
 \label{conditional}
G(\y)\,:=\,\mathrm{E}_\mu [\,h\,\vert \,\y]\,=\,\int h(\bl \y+\bs\,\z\,)\,\pi(d\z\vert\y)\,,
\end{equation}
where $\pi(d\z\vert\y)$ is the conditional probability on $\z$ given $\y$. They propose to evaluate the integral \eqref{conditional} by Monte-Carlo sampling on $\z$. But this sample \emph{depends} on $\y$ and therefore a sample has to be generated for each $\y$.

Our third contribution and novelty is  to exploit that if $h$ is a polynomial and $\mu$ is the uniform distribution on 
$\E_n$, then after the simple scaling $\v\to\z/\sqrt{1-\Vert\y\Vert^2}$, $\pi(d\v\vert\y)$ in \eqref{conditional} is the uniform distribution on $\E_{n-m}$. Therefore as the integrand is a polynomial in $\v$ of fixed degree, $G(\y)$ is a polynomial in the
$m+1$ variables $\y$ and $\sqrt{1-\Vert\y\Vert^2}$.
Its coefficients can be obtained exactly, e.g. by direct integration term by term after expansion of the integrand in the monomial basis.
 Alternatively,  $G$ can be expressed directly in terms of $h$ via a cubature formula on $\E_{n-m}$. Importantly, the cubature \emph{does not depend} on $\y$. Then for optimization on say
 $\bS^{n-1}$ or $\E_n$, instead of minimizing $h$, one proposes to minimize $G(X)=G(\bl^T\x)$ on  $\E_m$. This in turn is equivalent to minimizing a related function $\hat{f}(X,\vert Y\vert)$ on  $(X,Y)\in\bS^m$ for some polynomial $\hat{f}\in\R[X_1,\ldots X_m,Y]$.

%\texttt{ifacconf.tex}. You will also need the class file
%\texttt{ifacconf.cls}. Both files are available on the IFAC web site.
%
%Please stick to the format defined by the \texttt{ifacconf} class, and
%do not change the margins or the general layout of the paper. It
%is especially important that you do not put any running header/footer
%or page number in the submitted paper.\footnote{
%This is the default for the provided class file.}
%Use \emph{italics} for emphasis; do not underline.
%
%Page limits may vary from conference to conference. Please observe the 
%page limits of the event for which your paper is intended.
%

\section{Exploiting sparsity for optimization}
\label{exact-sparsity}
\subsection{Notation and definitions}

Let $C^1(\R^n)$ be the space of continuously differentiable functions on $\R^n$.
For any two vector $\x,\y\in\R^n$ denote by $\x\cdot\y$ their usual scalar product.
Given a vector space $V\subset\R^n$ denote by $V^\perp$ its orthogonal complement, i.e.,
$V^\perp=\{\,\y \in\R^n:\x\cdot\y=0\,,\:\forall \x\in V\,\}$.

The following result is relatively straightforward and its proof is omitted.
 \begin{prop}
 \label{prop1}
 Let $f:\R^n\times\R^m:\to\R$, $(\x,\y)\mapsto f(\x,\y)$, be continuously differentiable, and assume that 
 $\nabla_{\y} f(\x,\y)=0$ for all $(\x,\y)$.  Then with $\y_0\in\R^m$ fixed, arbitrary:
 \begin{equation}
 f(\x,\y)\,=\,f(\x,\y_0)\,=:\, g(\x)\,,\quad \forall (\x,\y)\in\R^n\times\R^m\,,
 \end{equation}
 and $g:\R^n\to\R$ is  continuously differentiable.
   \end{prop}
  
\subsection{Exploiting sparsity}

Let $h\in\R[\x]=\R[x_1,\ldots,x_n]$ and let 
\[\bl\,:=\,\left[\bl_1,\bl_2,\ldots,\bl_m\right]\,\in\,\R^{n\times m}\,,\]
for some $m$ linearly independent column vectors $\bl_1,\ldots\bl_m\in\R^n$.
Let $\x\cdot\y$ denote 
the usual scalar product of $\x,\y\in\R^n$.
\begin{thm}
\label{th1}
Let $h\in C^1(\R^n)$. Then the two statements below are equivalent: 

(a) There exists $f\in C^1(\R^m)$ and $(\bl_i)_{i=1,\ldots,m}\subset\R^n$ such that 
$h(\x)=f(\bl^T\,\x)=f(\bl_1\cdot\x,\ldots,\bl_m\cdot\x)$ for all $\x\in\R^n$.

(b) There exists an $m$-dimensional vector space $V\subset\R^n$ such that $\nabla h(\x)\in V$ for 
all $\x\in\R^n$.
\end{thm}
\begin{pf}
 (a) $\Rightarrow$ (b) is straightforward as
 \[\nabla h(\x)\,=\,\bl\, \nabla f(\bl^T\,\x)\,=\,\sum_{i=1}^m \bl_i\,\frac{\partial f(X)}{\partial X_i},\quad\forall \x\in\R^n\,.\]
 where $X_i=\bl_i\cdot\x$, $i=1,\ldots,m$.
 %, and $\bl\in\R^{m\times n}$ is the matrix with row vectors $\bl_i\in\R^n$, $i=1,\ldots,m$.
 Equivalently, $\nabla h(\x)\in V$ for all $\x\in\R^n$, where $V:={\rm Span}(\bl_1,\ldots,\bl_m)\subset\R^n$,
 which is clearly statement (b).

 (b) $\Rightarrow$ (a). Let $V\subset\R^n$ have dimension $m<n$ and let $(\bl_i)_{i=1,\ldots,m}$ be a basis of $V$. Similarly, let $(\bs_j)_{j=1,\ldots n-m}\subset\R^n$ be a basis of $V^\perp$ and write
$\x =\bl\,\u+\bs\,\v$,  with matrices $\bl=[\bl_1,\ldots\,\bl_m]\in\R^{n\times m}$ and $\bs=[\bs_1,\ldots\,\bs_{n-m}]\in\R^{n\times (n-m)}$, and such that $\bs_i^T\bl=0$ for all $i=1,\ldots,n-m$.
  Notice that
 \begin{equation}
 \label{ortho}
 \u\,=\,(\bl^T\,\bl)^{-1}\bl^T\,\x\,;\quad \v\,=\,(\bs^T\,\bs)^{-1}\bs^T\,\x\,.
 \end{equation}
 Hence write $h(\x)$ as
 \[h(\bl\,\u+\bs\,\v)\,=:\,\phi(\u\,,\v)\,=\,\phi((\bl^T\,\bl)^{-1}\bl^T\,\x\,,\,(\bs^T\,\bs)^{-1}\bs^T\,\x\,),\]
 for some function $\phi:\R^n\to\R$. Then  $\phi\in C^1(\R^n)$ follows from
 $h\in C^1(\R^n)$. Next, by the chain rule of differentiation:
 \[\nabla h(\x)\,=\,\bl\,(\bl^T\,\bl)^{-1}\nabla_{\u}\phi(\u,\v)+\bs\,(\bs^T\,\bs)^{-1}\nabla_{\v}\phi(\u,\v)\,.\]
 Observe that $\bs^T\cdot\nabla h(\x)=0$ for all $\x\in\R^n$ and all $i=1,\ldots,n-m$, because $\nabla h(\x)\in V$ for all $\x\in\R^n$. Hence 
 \[0\,=\,\bs^T\cdot\nabla h(\x)\,=\,\nabla_{\v}\phi(\u,\v)\,,\quad\forall \x\,\in\R^n\,,\]
 and therefore $\nabla_{\v}\phi(\u,\v)=0$, for all $\v\,\in\R^n$.
 By Proposition \ref{prop1} applied to $\phi$, $\phi(\u,\v)=\phi(\u,\v_0)$ for all $\u,\v$, where
  where $\v_0$ is arbitrary. Letting $\v_0:=0$ yields 
  \begin{eqnarray*}
  h(\x)\,=\,\phi(\u,\v)\,=\,\phi(\u,0)&=&\phi((\bl^T\bl)^{-1}\bl^T\x,0)\\
  &=&f(\tilde{\bl}_1\cdot\x,\ldots,\tilde{\bl}_m\cdot\x)\,,
  \end{eqnarray*}
  where $\tilde{\bl}_i\in\R^n$ is the $i$-th row of $(\bl^T\bl)^{-1}\bl^T$, $i=1,\ldots,m$, and 
  $f(X_1,\ldots,X_m):=\phi(X_1,\ldots,X_m,0)$ for all $X\in\R^m$.
 \end{pf}
\subsection{Detection of the sparse form}
 
 In this section we suppose that $h\in\R[\x]$ is a sparse polynomial but \emph{not } given
 in its sparse form  $\x\mapsto f(\bl^T\x)$ for some real matrix $\bl\in\R^{n\times m}$. In view of 
 Theorem \ref{th1}, it suffices to determine a basis of the $m$-dimensional subspace $V\subset\R^n$ 
 to  which $\nabla h$ belongs. 
 
 We consider a methodology inspired from Constantine et al.  \cite{active-set} but we here exploit that $h$ is a polynomial.
 Introduce a probability measure $\mu$ on a certain domain, e.g. the uniform  distribution  on $\mathcal{E}_n$: 
 
 (i) A first possibility consists in computing the $n\times n$ real symmetric matrix
 \[\M_\mu\,:=\,\mathrm{E}_\mu[\,\nabla h(\x)\,\nabla h(\x)^T\,]\,,\]
 and compute its SVD decomposition. In \cite{active-set} the function $h$ is \emph{not} a polynomial and therefore
 $\mathrm{E}_\mu$ must be approximated. Moreover $h$ is not necessarily sparse and the authors are interested 
 in approximating $h$ in the subspace $V$ generated by the eigenvectors associated  with the  largest eigenvalues of $\M_\mu$. In our setting, $\M_\mu$
 can be computed \emph{exactly} as one knows how to integrate \emph{exactly} a polynomial on $\mathcal{E}_n$, and
 $V$ is spanned by the eigenvectors associated with the zero-eigenvalues of $\M_\mu$.
 
 (ii) Another possibility is to consider a sample of $(m+1)$ i.i.d. vectors  $(\nabla h(\x(i))_{i\leq m+1}\subset\mathcal{E}_n$ randomly generated according to $\mu$, 
 and construct the empirical matrix 
 \begin{equation}
 \label{empirical}
 \bH_{m+1}:=[\nabla h(\x(1)),\cdots,\nabla h(\x(m+1))]\in\R^{n\times (m+1)}\,\end{equation}
 until one observes that ${\rm rank}(\bH_\ell^T\,\bH_\ell)=m$, $\ell=m,m+1$.
 \begin{thm}
\label{th2}
  Let 
 $\bH_k:=[\nabla h(\x(1)),\cdots,\nabla h(\x(k))]\in\R^{n\times k}$ be as in \eqref{empirical},
  and let
 $V:={\rm span}\{\nabla h(\x): \x\in \E_n\,\}$.
  Then ${\rm dim}(V)=m$ if and only, with probability $1$:
  \begin{equation}
 \label{prob-1}
{\rm rank}(\bH_\ell^T\,\bH_\ell)\,=\,m\,,\quad \ell=m,m+1\,.%\quad\mbox{and}\quad {\rm rank}(\bH_{m+1}\,\bH_{m+1}^T)\,=\,m\,.
\end{equation}
\end{thm}
 
\begin{pf}
The \emph{Only if part} is straightforward. Indeed in view of the definition of $V$,
suppose that ${\rm dim}(V)=m$,  and let
 $V^\perp$ denote its direct complement (hence of  dimension $n-m$). Then $\u^T\nabla h(\x(i))=0$ for all $\u\in V^\perp$ and all $i=1,\ldots,m$,
% Hence $\u^T\bH_m\bH_m^T=0$ for all $\u\in V^\perp$, 
 which implies
$ {\rm rank}(\bH_\ell^T\,\bH_\ell)\,\leq\,m$, $\ell=m,m+1$.

Next, observe that 
 \[(\bH^T_m\,\bH_m)_{i,j}\,=\,\nabla h(\x(i))^T\nabla h(\x(j))\,,\quad i,j\le m\,,\]
 and therefore 
 \begin{equation}
 \label{def-pm}
 {\rm det}(\bH^T_m\,\bH_m)\,=\,p_m(\x(1),\x(2),\ldots,\x(m))\,,\end{equation}
 for some polynomial $p_m\in\R[\u_1,\ldots,\u_m]$. As ${\rm dim}(V)=m$ then necessarily
 $p_m\neq0$.   
 Next, let $\mu^{\otimes m}$ be the product measure
 $\underbrace{\mu\otimes\mu\cdots\otimes\mu}_{m\:times}$ on $(\E_n)^m$.
 As $p_m\neq0$ then $\mu^{\otimes m}(\{\u: p_m(\u_1,\ldots,\u_m)=0\})=0$,
 or equivalently,
 with probability $1$, $p_m(\u_1,\ldots,\u_m)\neq0$, i.e., ${\rm det}(\bH^T_m\,\bH_m)\,\neq\,0$, and so
 ${\rm rank}(\bH_m^T\bH_m)=m$.
 %\[\mbox{${\rm rank}(\bH_m\bH_m^T)\,<\,m$ with probability $1$}\quad\Leftrightarrow  \quad 
 Next consider the case $\ell =m+1$. As ${\rm dim}(V)=m$ then necessarily the family
 $(\nabla h(\x(i)))_{i\leq m+1}$ is \emph{not} linearly independent and therefore
 ${\rm rank}(\bH_{m+1}^T\bH_{m+1})<m+1$, which from  what precedes,
 yields ${\rm rank}(\bH_{m+1}^T\bH_{m+1})=m$ with probability $1$.
 
 \emph{If part}. As above, let
 \begin{equation}
 \label{def-pk}
  {\rm det}(\bH^T_k\,\bH_k)\,=:\,p_k(\x(1),\x(2),\ldots,\x(k))\,,\quad k\in\N\,,\end{equation}
 for some polynomial $p_k\in\R[\u_1,\ldots,\u_k]$. 
 The condition
\[ \mbox{``with probability $1$,}\quad {\rm rank}(\bH_\ell^T\,\bH_\ell)\,=\,m\,,\quad \ell=m,m+1\,"\,,\]
 is equivalent to
 %\begin{equation} \label{condi}
  \[\mbox{``with probability $1$:}\quad \left\{\begin{array}{l}{\rm det}(\bH_{m}^T\,\bH_{m})\,>\,0\,,\mbox{ and}\\
 {\rm det}(\bH_{m+1}^T\,\bH_{m+1})\,=\,0\,,\end{array}\right."\]
 %\end{equation}
 which in turn is equivalent to 
  \begin{equation}
 \label{condi}
 p_m\,\neq\,0\mbox{ and }p_{m+1}\,=\,0\,,\end{equation}
 with $p_m$ as in \eqref{def-pk}.
 %To emphasize the dependence on the sample we write
%\[\bH_m=\bH_m(\x(1),\ldots,\x(m))\,;\quad \bH_{m+1}=\bH_m(\x(1),\ldots,\x(m),\x(m+1))\,.\]
 %\quad\mbox{and}\quad {\rm rank}(\bH_{m+1}\,\bH_{m+1}^T)\,=\,m\,.
 %As $\lambda^{\otimes m+1}$ is the Lebesgue measure on $([0,1]^n)^{m+1}$,
 The condition %on the right in \eqref{condi} implies 
 $p_{m+1}=0$, i.e., 
 \[%\begin{eqnarray*}
 {\rm det}(\bH_{m+1}(\u_1,\ldots,\u_{m+1})^T\bH_{m+1}(\u_1,\ldots,\u_{m+1}))\,=\,0\,,\]
 for all $\u:=(\u_1,\ldots\u_{m+1})\in (\E_n)^{m+1}$, implies
 that  there exists  a vector $0\neq q^{\u}\in\R^{m+1}$ such that
 \[\bH_{m+1}(\u_1,\ldots,\u_{m+1})\,q^{\u}\,=\,0\,,\quad\forall \u\in (\E_n)^{m+1}\,.\]
 
 %$ there exists \]
 That is,
\[ \sum_{i=1}^{m+1}q^{\u}_i\,\nabla h(\u_i)\,=\,0\,,\quad\forall \u\in (\E_n)^{m+1}\,.\]
Next, let $S:=\{\u_{m+1}\in \E_n: q^\u_{m+1}=0\}$ and $\Theta=(\E_n)^{m}\times S$, so that
$\mu^{\otimes (m+1)}(\Theta)=\mu(S)$. Hence
\[ \sum_{i=1}^{m}q^{\u}_i\,\nabla h(\u_i)\,=\,0\,,\quad\mbox{for all $\u\in \Theta$}\,.\]
 %\[\mbox{${\rm rank}(\bH_m\bH_m^T)\,<\,m$ with probability $1$}\quad\Leftrightarrow  \quad p\,=\,0\,.\]
Next, from ${\rm det}(\bH_{m}(\u_1,\ldots,\u_m)^T\,\bH_{m}(\u_1,\ldots,\u_m))>0$,
we also deduce that 
\begin{equation}
\label{aux}
 \sum_{i=1}^{m}q^{\u}_i\,\nabla h(\u_i)\,\neq0\,,\quad\mbox{for a.a. $\u\in (\E_n)^{m+1}$}\,.\end{equation}
This yields  $0\,=\,\mu^{\otimes (m+1)}(\Theta)=\mu(S)$.
Therefore, letting $\u(\x):=(\u_1,\ldots,\u_{m},\x)\in (\E_n)^m\times (\E_n\setminus S)$, %for all $\x\in [0,1]^n\setminus S$,
\[\nabla h(\x)\,=\,\frac{1}{q^{\u(\x)}_{1}}\,
\sum_{i=1}^{m} q^{\u(\x)}_i\,\nabla h(\u_i)\,,\]
for all $\x\in \E_n\setminus S$, and  all $(\u_1,\ldots,\u_{m})\in (\E_n)^m$.
Hence with $(\u_2,\ldots,\u_{m+1})\in (\E_n)^m$, fixed, arbitrary:
\begin{equation}
\label{fin-2}
\nabla h(\x)\in\,{\rm span}(\nabla h(\u_1),\ldots,\nabla h(\u_{m}))\,=:\,V\,,
\end{equation}
for all $\x\in \E_n\setminus S$,
and $V$ is an $m$-dimensional vector space. To show that \eqref{fin-2} holds for all $\x\in\R^n$, observe that
\begin{equation}
 \label{fin}
 \v^T\nabla h(\x)\,=\,0\,,\quad\forall \x\in \E_n\setminus S\,,\:\forall \v\in V^\perp.
 \end{equation}
Hence for fixed $\v\in V^\perp$, the polynomial $\x\mapsto \v^T\nabla h(\x)$ vanishes on $\E_n\setminus S$
with $\mu(S)=0$, which implies that $\v^T\nabla h(\x)$ vanishes on the whole $\E_n$ and hence 
on the whole $\R^n$. As this is true for an arbitrary $\v\in V^\perp$, we obtain that $\nabla h(\x)\in V$ for all $\x\in\R^n$.
\end{pf}

In practice,  Theorem \ref{th2} is used as follows:
\begin{itemize}
\item Samples $k$ points $(\x(i))_{i\leq k}$ according to $\lambda$ on $[0,1]^n$.
\item Do the SVD decomposition of the real symmetric matrices $\bH_{k-1}^T\bH_{k-1}$.
and $\bH_{k}^T\bH_{k}$.
\item If ${\rm rank}(\bH_k^T\bH_k)\,\neq\,{\rm rank}(\bH_{k-1}^T\bH_{k-1})$ then set $k:=k+1$ and repeat.
\item If ${\rm rank}(\bH_k^T\bH_k)\,=\,{\rm rank}(\bH_{k-1}^T\bH_{k-1})$ then stop.\\
Set
$V:={\rm span}\{\nabla h(\x(1),\ldots,\nabla h(\x(k-1))\}$.
\end{itemize}
 \subsection{Some applications in Optimization}
 
 \subsubsection{Optimization on the Euclidean unit sphere}
 A first application was developed in \cite{lass-sphere} for optimization on the Euclidean unit sphere $\mathbb{S}^{n-1}$.
 Namely, let $h,f$ and $\bl$ be as in Theorem \ref{th1}(a). Then it was shown in \cite{} that
 \begin{eqnarray}
 \label{a1}
 \rho&=&\min_{\x}\{\,h(\x)\,:\: \x\in\mathbb{S}^{n-1}\,\}\\
 \label{a2}
 &=&\min_{\y}\{\,f(\mathbf{L}_1\cdot \y_1,\ldots,\mathbf{L}_m\cdot \y_m)\,:\:\y\in\mathcal{E}_m\,\}\,,\
 \end{eqnarray}
 with $\mathbf{L}_i$ is the $ith$-row of the matrix $(\bl^T\bl)^{1/2}$, $i=1,\ldots,m$. In fact all points $\x^*\in\mathbb{S}^{n-1}$ that satisfy 
 the standard second-order necessary conditions of optimality for problem \eqref{a1} are in one-to-one correspondence
 with the points $\y^*\in\mathcal{E}_m$ that satisfy the standard second-order necessary conditions of optimality for 
 problem \eqref{a2}. 
 
 Hence in this case one has replaced optimization of the $n$-variate  polynomial $h$ on the non convex set $\mathbb{S}^{n-1}$ by optimization of the $m$-variate polynomial $f$ of same degree on the (convex) unit Euclidean ball. If $m\ll n$
 then it yields drastic computational savings.
 
 \subsubsection{Optimization on a polytope}
 Next, let $\bom=\{\,\x\in\R^n_+:\:\A\x = \b\,\}$
 for some real matrix $\A\in\R^{s\times n}$, and  consider the optimization problem:
 \begin{equation}
 \label{polytope}
 \rho\,=\,\min_{\x}\{\,h(\x)\,:\: \x\in\bom\,\}\,.
 \end{equation}
 \begin{thm}
 Let $h$ and $f$ be as in Theorem \ref{th1}(a) with $\bl\in\R^{n\times m}$, and let
 $(\bla_i,\u_i)_{i\in I}$ be a set of generators of the polyhedral convex cone
 \begin{equation}
 \label{cone-C}
 C\,:=\,\{(\bla,\u)\in\R^{s}\times \R^m: \:\A^T\bla\,\geq\,\bl\,\u\,\}\,.\end{equation}
 Then with $\rho$ as in \eqref{polytope}
 \begin{equation}
 \rho \,=\,\min_{X\,\in\,\R^m}\,\{\,f(X)\::\: \u_i\cdot X\,\leq\,\bla_i\cdot\b\,,\quad \forall\,i\in I\,\}
\end{equation}
 \end{thm}
 
\begin{pf}
 Let $X$ be fixed. By Farkas Lemma,  
 \[\emptyset\,\neq\,\{\x:\: \bl^T\,\x\,=\,X\,;\: \x\in\bom\}\,\Leftrightarrow\,\u\cdot X\,\leq\,\bla\cdot\b\,,\]
 for all $(\bla,\u)\,\in\,C$,
 which in turn is equivalent to $\u_i\cdot X\,\leq\,\bla_i\cdot\b$ for all $(\bla_i,\u_i)_{i\in I}$.
 Then observe that
 \begin{eqnarray}
 \nonumber
 \P:\quad \rho&=&\min_{\x}\{\,h(\x)\,:\: \x\in\bom\,\}\\
 \label{optim-P}
 &=& \min_{\x,X}\,\{\,f(X)\,:\: X\,=\,\bl^T\x\,;\: \x\in\bom\,\}\\
 \nonumber
 &=&\min_{X}\,\{\,f(X)\::\: \u_i\cdot X\,\leq\,\bla_i\cdot\b\,,\quad \forall\,i\in I\,\}\,.
 \end{eqnarray}
  \end{pf}
 Notice that one has replaced an $n$-dimensional optimization problem on the polyhedron $\bom\subset\R^n$ 
 by an $m$-dimensional optimization problem on the polyhedron 
 $\bom_m:=\{X\in\R^m: \u_i\cdot X\leq\bla_i\cdot\b\,,\: i\in I\,\}\subset\R^m$. 
 Of course this transformation requires to compute as a pre-requisite step, all generators of the convex cone $C$ in \eqref{cone-C}.
 If one wants to avoid this, one possibility is to proceed as follows:
 
 $\bullet$ Start with a set $I_0:=\{(\bla_0,\u_0)\}$ for some $(\bla_0,\u_0)\in C$, and set $k=0$.
 
 $\bullet$ Step $k$. Solve 
 \[\mathbf{P}_k:\quad \tau_k=\min_X\,\{f(X)\::\: \u_i\cdot X\leq\bla_i\cdot\b\,,\quad i\in I_k\,\}\,,\]
 to obtain $X^*_k\in\R^m$. Next, solve the linear program 
 %\[\mathbf{L}_k:\quad 
 \[\begin{array}{rl}
 \tau=\displaystyle\min_{\bla^+,\bla^-,\u^+,\u^-}\,&\{(\bla^+-\bla^-)\cdot\b-(\u^+-\u^-)\cdot X^*_k\,:\\
 &((\bla^+-\bla^-),(\u^+-\u^-))\in C\,;\\
 &\sum_t\lambda^+_t+\lambda^-_t+\sum_j(u^+_j+u^-_j)=1\,\}\,.\end{array}\]
 If $\tau=0$ then stop.
 Otherwise set $I_{k+1}:=I_k\cup\{(\bla_*,\u_*)\}$ for an optimal solution
 $(\bla^+_*-\bla^-_*,\u^+_*-\u^-_*)$, set $k:=k+1$  and go to step $k$.
 
 With this strategy one has to solve a sequence of optimization problems $(\mathbf{P}_k)_{k\in\N}$
 with same criterion $f$, but on tighter and tighter outer approximations of
 the convex polyhedron $\{X\in\R^m: \u_i\cdot X\leq \bla_i\cdot\b,\,,\: i\in I\}$. So the overall complexity
 of this algorithm is governed by the computational complexity of problem $\mathbf{P}_k$.
 
 For simple sets $\bom$ like the canonical simplex or the unit box, the cone $C$ has a simple expression.
 \subsubsection{On the canonical simplex $\bom=\{\x\in\R^n_+:\e\cdot\x =1\}$.}
 $C\,=\,\{(\lambda,\u): \lambda\, \e\geq\,\bl\,\u\}$.
 \subsubsection{On the Box $\bom=[-1,1]^n$.} %Another important special case is the canonical Box
 % $\bom:=\{\x\in\R^n: \vert x_i\vert \leq 1\,;\:i\in [n]\}$.
 $C=\{(\bla^+,\bla^-\geq0,\u): \bla^+-\bla^-=\bl\,\u\,\}$ and so $(\bla^++\bla^-)\cdot\e\,=\Vert\bl \u\Vert_1$.
 %,  and therefore $\P$ is equivalent to solving
  %\[\min_{X}\,\{\,f(X)\::\: \u^T\,X\,\leq\,\Vert\bl\,\u\,\Vert_1\,,\:\forall \u\in\mathbb{S}^{m-1}\}\,.\]
 
 \section{Approximate sparsity}
 
 In this section $h\in\R[\x]$ and we now assume that $h$ is not exactly in the form $f(\bl^T\,\x)$ for some $\bl\in\R^{n\times m}$. Let %$\E_n:=[-1,1]^n$, 
 $\mu_n$ be the uniform distribution on $\E_n$ and let
 \[\M(\mu)\,:=\,\mathrm{E}_{\mu_n}[\,\nabla h(\x)^T\,\nabla h(\x)\,]\,=\,[\bl,\bs]\,\left[\begin{array}{cc}\Lambda_1 &0\\ 0 &\Lambda_2\end{array}\right]\,[\bl,\bs]^T\]
 where now $\bl=[\bl_1,\ldots,\bl_m]\in\R^{n\times m}$ (resp. $\bs=[\bs_1,\ldots,\bs_{n-m}]\in\R^{n\times (n-m)}$ is the matrix eigenvectors of $\M(\mu_n)$ associated with
 the first $m$ (nonnegative) eigenvalues $\lambda_1,\ldots,\lambda_m$ 
 (resp. the remaining $n-m$ eigenvalues $\lambda_{m+1},\ldots,\lambda_n$) arranged in decreasing order and which also form the diagonal elements of the diagonal matrices
 $\Lambda_1$ and $\Lambda_2$, respectively. The vectors
 $\bl_j,\bs_j$ form an orthonormal basis. Therefore  if one writes $\x=\bl\, \y+\bs\, \z$ with $\y\in\R^m$ and  $\z\in\R^{n-m}$, then
  \[\Vert \x\Vert^2\,=\,\Vert \y\Vert^2+\Vert \z\Vert^2\,,\]
 and so the support of the marginal $\pi_\y$ of $\mu_n$ is $\E_m$,
 with density (w.r.t. Lebesgue) the pushforward of $\mu_n$ by its projection on $\E_m$.
 The conditional $\pi(d\z\vert \y)$ is the uniform probability distribution on 
 the  ball $\E_{n-m}(\y):=\{\z:\Vert \z\Vert^2\leq 1-\Vert \y\Vert^2\}$.
 Proceeding as in \cite{active-set}, introduce the function
  %\footnote{$f$ is not a polynomial because the sample $(Z(j))_{j\leq N}$ depends on $X$.} 
  $f:\R^m\to \R$, defined by
 \begin{eqnarray}
 \label{cond}
 f(\y)&:=&\mathrm{E}[h\vert \y]\\
 \nonumber
 &=&\displaystyle\int_{\E_{n-m}(\y)} h(\bl \,\y+\bs\,\z)\,\pi(d\z\vert\y)\,,\quad\forall\y\in\E_m\,.
 \end{eqnarray}
 Then the idea promoted in \cite{active-set} for some 	applications, is to approximate $h$
 on $\E_n$ with the function $\hat{h}(\x):=f(\bl^T\x)$. The rationale being:
 \begin{thm}(\cite[Theorem 3.1]{active-set})%\cite{active-set}[Theorem 3.1])
 \label{th-active-set}
 With $\mu_n$ being the uniform probability distribution on $\E_n$, and $\hat{h}(\x)=f(\bl^T\x)$, with $f$ as in 
 \eqref{cond},
 \begin{equation}
 \label{small}
 \mathrm{E}_{\mu_n}[\,(h-\hat{h})^2]\,\leq\,C\,(\lambda_{m+1}+\ldots,+\lambda_n)\,.
\end{equation}
where the constant $C$ does not depend on $h$.
\end{thm}
So in view of \eqref{small}, if the remaining eigenvalues
$\lambda_{m+1},\ldots,\lambda_n$ are small then 
$\hat{h}$ provides a good approximation  of $h$ 
in $L^2(\E_n)$.
 \subsection*{Exact computation of the approximand $f$}
 Observe that $\pi(d\z\vert\y)=d\z/C_{n-m}(1-\Vert\y\Vert^2)^{(n-m)/2}$ on $\E_{n-m}(\y)$, for a constant $C_{n-m}$.
 Therefore by  doing the change of variable $\v:=\z/\sqrt{1-\Vert\y\Vert^2}\in\E_{n-m}$,  and letting $\tau(\y):=1-\Vert\y\Vert^2$,
 \eqref{cond} reads:
 \begin{eqnarray}
 \nonumber
  \mathrm{E}[h\vert \y]&=&\frac{1}{C_{n-m}}\,
 \displaystyle\int_{\E_{n-m}} 
 h(\bl \,\y+\tau(\y)^{1/2}\,\bs\,\v\,) \,d\v\,,\\
 \label{E[h-y]}
 &=&\displaystyle\int_{\E_{n-m}}  h(\bl \,\y+\tau(\y)^{1/2}\,\bs\,\v\,) \,d\mu_{n-m}(\v)\,,
 \end{eqnarray}
 for all $\y\in\E_m$. 
 Observe that the integrand $\v\mapsto h(\bl \,\y+\tau(\y)^{1/2}\,\bs\,\v)$ is a polynomial of fixed degree, say $d$, in $\v$. Therefore it can be integrated exactly on $\E_{n-m}$. Equivalently one can
 also use a degree-$d$ cubature rule for Lebesgue measure on $\E_{n-m}$ to obtain:
 \begin{equation}
 \label{cubature}
 f(\y)\,=\,\sum_{j=1}^r
 \theta_j\,h(\bl \,\y+\tau(\y)^{1/2}\bs\,\v_j\,)\,,
 \end{equation}
  for some positive weights $(\theta_j)$ and cubature points
  $(\v_j)\subset\E_{n-m}$. Importantly, and 
   in contrast to the function $G(\y)$ in \cite[(3.10)]{active-set},
   the cubature points $(\v_j)$ do \emph{not} depend on $\y$ and so can be computed once and for all\footnote{In \cite{active-set} the integral $\mathrm{E}[h\vert \y]$
   has to be computed via Monte-carlo sampling with a different sample for each $\y$.}. Notice that
   $f$ is a polynomial in the $(m+1)$ variables 
   $(y_1,\ldots,y_m,\sqrt{1-\Vert\y\Vert^2})$, i.e., 
   $f\in\R[\y,\sqrt{1-\Vert\y\Vert^2}]$. Next, again following \cite{active-set} we approximate $h$ on $\E_n$ with
   $h(\x)\approx\hat{h}(\x)\,:=\,f(\bl^T\x)$, i.e.:
   \begin{equation}
  \label{approx-h}
 \hat{h}(\x)\,=\,\sum_{j=1}^r
 \theta_j\,h(\bl \,\bl^T\x+\tau(\bl ^T\x)^{1/2}\bs\,\v_j\,)\,.
  \end{equation}
  Hence letting $X:=\bl^T\x$ and using the orthogonality of the $(\bl_j)$,  we obtain  $X\in\E_m$, and 
 \begin{equation}
 \label{approx-2}
 f(X)\,=\,\sum_{j=1}^r
 \theta_j\,h(\bl \,X+(1-\Vert X\Vert^2)^{1/2}\,\bs\,\v_j\,)\,.
 \end{equation}
  Next, introduce  the polynomial $\hat{f}\in\R[X,Y]$ with
 \begin{equation}
 \label{approx-3}
 \hat{f}(X,Y)\,:=\,\sum_{j=1}^r
 \theta_j\,h(\bl \,X+Y\,\bs\,\v_j\,)\,,
 \end{equation}
 for all $(X,Y)\in\R^{m+1}$, and let $Y^2=1-\Vert X\Vert^2$ so that $(X,Y)\in\bS^m$ whenever $X\in\E_m$.
 Hence whenever $\x\in\E_n$, then $(X,Y)\in\bS^m$, and
 \begin{equation}
 \label{final-form}
 \hat{h}(\x)\,=\,\hat{f}(\bl^T\x,\sqrt{1-\Vert\bl^T\x\Vert^2})\,=\,\hat{f}(X,\vert Y\vert)\,\:\mbox{on $\bS^m$}\,.
 \end{equation}
 
\subsection*{Approximate sparse optimization on $\E_n$ or $\bS^{n-1}$}So when the $n-m$ remaining eigenvalues $(\lambda_{m+1},\ldots,\lambda_n)$ are small compared to the first $m$ ones, Theorem \ref{th-active-set} suggests to consider 
replacing $h$ with $\hat{h}$ in the initial optimization problem $\P$. As we next show, when $\bom=\bS^{n-1}$ or $\bom=\E_n$,
the resulting problem is equivalent to solving:
\begin{equation}
\label{final-Q}
 \Q:\quad\rho=\min_{(X,Y)}\{\,\hat{f}(X,\vert Y\vert):\: (X,Y)\,\in\, \bS^m\,\}\,,
\end{equation}
an $(m+1)$-variables optimization problem. Note that $\hat{f}(X,\vert Y\vert)$ is not a polynomial but 
$\rho=\min [\,\rho^+,\rho^-\,]$ with
\begin{eqnarray*}
\rho^+&=&\min_{X,Y}\{\,\hat{f}(X,Y):\: (X,Y)\,\in\, \bS^m\,;\,Y\geq0\}\\
\rho^-&=&\min_{X,Y}\{\,\hat{f}(X,-Y):\: (X,Y)\,\in\, \bS^m\,;\,Y\leq0\}\,.
\end{eqnarray*}
So to solve $\Q$ and obtain $\rho$, one has to solve two \emph{polynomial} optimization problems of same type as $\P$ but on $\bS^m$, hence of much lower dimension when $m\ll n$.

\begin{lem}
 Let $\hat{h}$ be as in \eqref{approx-h},  $\hat{f}$ as in \eqref{approx-3}, and let $\rho=\min[\,\rho^+,\rho^-]$. Then 
 \begin{equation}
 \label{equality}
 \min \{\,\hat{h}(\x): \x\in\bS^{n-1}\}\,=\,\min \{\,\hat{h}(\x): \x\in\E_n\}\,=\,\rho\end{equation}
 \end{lem}
\begin{pf}
Let $\tau:=\min \{\,\hat{h}(\x): \x\in\E_n\}$ and
let $\x^*:=\arg\min \{\,\hat{h}(\x): \x\in\E_n\}$ so that
$\hat{h}(\x^*)= \tau$. Write $\x^*=\bl\,\y^*+\bs\,\z^*$
so that $\Vert\x^*\Vert^2=\Vert \y^*\Vert^2+\Vert\z^*\Vert^2\leq1$.
Next, let $\tilde{\x}:=\bl\,\y^*+r\cdot\bs\,\z^*$ so that
$\Vert\tilde{\x}\Vert^2=\Vert \y^*\Vert^2+r^2\,\Vert\z^*\Vert^2$,
and choose $r$ such that $\tilde{\x}\in\bS^{n-1}$. Then
$\bl^T\tilde{\x}=\bl^T\x^*$ and therefore $\hat{h}(\tilde{\x})=\hat{h}(\x^*)=\tau$, which yields the 
first equality in \eqref{equality}. 
It remains to prove that $\rho=\tau$.

It is clear that $\tau\geq\rho$ as 
$(X,Y):=(\bl^T\x,\sqrt{1-\Vert X\Vert^2})\in\bS^m$ and
$\hat{f}(X,Y)=\hat{h}(\bl^T\x)$ whenever $\x\in \E_n$.
For the converse, assume that $\rho=\rho^+$ with an optimal solution 
$(X^*,Y^*)\in\bS^m$ and $Y^*\geq0$. Let
$\x:=\bl\,X$ so that $\Vert\x\Vert=\Vert X\Vert$ as $\bl^T\bl=\mathrm{I}_m$. Hence $\x\in\E_n$.
Moreover $\bl^T\x=\bl^T\bl\,X=X$ and therefore by
\eqref{approx-h}-\eqref{approx-2}, $\hat{h}(\bl^T\x)=\hat{f}(X,Y)=\hat{f}(X,\vert Y\vert)$, which proves that $\tau\leq \rho^+$. The proof when $\rho=\rho^-$ being similar is omitted.
\end{pf}
Of course the rationale for solving $\Q$ instead of  $\P$ is based on 
Theorem \ref{th-active-set}, assuming that
$\sum_{j=m+1}^n\lambda_j\,(\mathrm{E}_{\mu_n}[\nabla h\nabla h^T])$ is small.
But the approximation in Theorem \ref{th-active-set} in only in $L^2(\E_n)$ and not in $L^\infty(\E_n)$ (or equivalently in the sup-norm). This is why we have not provided
an error analysis which remains to be done.

Notice that if $\lambda_j\,(\mathrm{E}_{\mu_n}[\nabla h\nabla h^T])=0$ for all $j>m$, then one retrieves the problem of Section \ref{exact-sparsity}. Indeed in $h(\x)=h(\bl \y+\bs\,\z)$ one has 
$\z=0$, and therefore in \eqref{cond} and \eqref{E[h-y]}, 
\[f(\y)\,=\,\mathrm{E}_{\mu_n}[h\vert\y]\,=\,h(\l\,\y)\,.\]
So for instance when $\bom=\bS^{n-1}$ and $h(\x)=f(\bl^T\x)$ for some $\bl\in\R^{n\times m}$,
the  sparse problem $\Q=\min \{f(\L\,\y):\y\in\E_m\}$ shown to be strictly equivalent to $\P$ in \cite{lass-sphere}, 
is the limit case of  $\Q$ in \eqref{final-Q} when 
$\sum_{j=m+1}^n\lambda_j(\mathrm{E}_{\mu_n}[\nabla h\nabla h^T])=0$.

\end{document}